\documentclass{article}
\usepackage{amsmath}
\usepackage{amsfonts}
\usepackage{amssymb}
\usepackage{graphicx}

\providecommand{\U}[1]{\protect\rule{.1in}{.1in}}
\input{tcilatex}

\begin{document}

\title{Projective Modules over Quantum Projective Line\thanks{%
This work is part of the project supported by the grant
H2020-MSCA-RISE-2015-691246-QUANTUM DYNAMICS.}}
\date{}
\author{Albert Jeu-Liang Sheu\thanks{%
The author would like to thank the Mathematics Institute of Academia Sinica
for the warm hospitality and support during his visit in the summer of 2015.}
\\
%EndAName
{\small Department of Mathematics, University of Kansas, Lawrence, KS 66045,
U. S. A.}\\
{\small e-mail: asheu@ku.edu}}
\maketitle

\begin{abstract}
Taking a groupoid C*-algebra approach to the study of the quantum complex
projective spaces $\mathbb{P}^{n}\left( \mathcal{T}\right) $ constructed
from the multipullback quantum spheres introduced by Hajac and
collaborators, we analyze the structure of the C*-algebra $C\left( \mathbb{P}%
^{1}\left( \mathcal{T}\right) \right) $ realized as a concrete groupoid
C*-algebra, and find its $K$-groups. Furthermore after a complete
classification of the unitary equivalence classes of projections or
equivalently the isomorphism classes of finitely generated projective
modules over the C*-algebra $C\left( \mathbb{P}^{1}\left( \mathcal{T}\right)
\right) $, we identify those quantum principal $U\left( 1\right) $-bundles
introduced by Hajac and collaborators among the projections classified.
\end{abstract}

\section{Introduction}

In the theory of noncommutative topology or geometry \cite{Conn}, a
generally noncommutative C*-algebra $\mathcal{A}$ is viewed as the algebra $%
C\left( X_{q}\right) $ of continuous functions on a virtual spatial object $%
X_{q}$, called a quantum space. Many interesting examples of quantum spaces
have been constructed with a topological or geometrical motivation, and
analyzed in comparison with their classical counterpart. Different
topological or geometrical viewpoints of the same object may give rise to
different quantum versions of quantum spaces. For example, quantum
odd-dimensional spheres and associated complex projective spaces have been
introduced and studied by Soibelman, Vaksman, Meyer, and others \cite%
{VaSo,Me} as $\mathbb{S}_{q}^{2n+1}$ and $\mathbb{C}P_{q}^{n}$ via quantum
universal enveloping algebra approach, and by Hajac and his collaborators
including Baum, Kaygun, Matthes, Pask, Sims, Szyma\'{n}ski, Zieli\'{n}ski,
and others \cite{BaHaMaSz,HaMaSz,HaKaZi,HaPaSiZi} as $\mathbb{S}_{H}^{2n+1}$
and $\mathbb{P}^{n}\left( \mathcal{T}\right) $ via a multi-pullback and
Toeplitz algebra approach.

Recall that the concept of a vector bundle $E$ over a compact space $X$ can
be reformulated in the noncommutative context as a finitely generated
projective left module $\Gamma\left( E_{q}\right) $ over $C\left(
X_{q}\right) $, viewed as the space of continuous cross-sections of some
virtual noncommutative or quantum vector bundle $E_{q}$ over $X_{q}$, as
suggested by Swan's work \cite{Swan}. Based on the strong connection
approach to quantum principal bundles \cite{Haja:sc} for compact quantum
groups \cite{Woro,Wo:cm}, Hajac and his collaborators introduced quantum
line bundles $L_{k}$ of degree $k$ over $\mathbb{P}^{n}\left( \mathcal{T}%
\right) $ as some rank-one projective modules realized as spectral subspaces 
$C\left( \mathbb{S}_{H}^{2n+1}\right) _{k}$ of $C\left( \mathbb{S}%
_{H}^{2n+1}\right) $ under a $U\left( 1\right) $-action, and analyzed them
via pairing of cyclic cohomology and K-theory \cite{HaMaSz,HaPaSiZi}. In
particular, it was found that $L_{k}$ is not stably free unless $k=0$,
revealing some information about the $K_{0}$-group of $C\left( \mathbb{P}%
^{n}\left( \mathcal{T}\right) \right) $. On the other hand, even for the
most crucial case of $n=1$, the $K_{0}$-group was not fully determined
despite the significant progress made in the 2003 paper \cite{HaMaSz}.

Going beyond the $K$-theoretic study of C*-algebras that classifies finitely
generated projective modules only up to stable isomorphism, some successes
have been achieved in the study of cancellation problem, made popular by
Rieffel \cite{Ri:dsr,Ri:ct}, that classifies finitely generated projective
modules up to isomorphism, for some quantum algebras \cite%
{Ri:ct,Ri:pm,Sh:ct,Bach,Pete}. It is of interest and a natural question to
classify finitely generated projective modules over $C\left( \mathbb{P}%
^{n}\left( \mathcal{T}\right) \right) $ and identify the line bundles $L_{k}$
among them, beside finding the $K$-groups of $C\left( \mathbb{P}^{n}\left( 
\mathcal{T}\right) \right) $.

In this paper, we use the powerful groupoid approach to C*-algebras
initiated by Renault \cite{Rena} and popularized by Curto, Muhly, and
Renault \cite{CuMu,MuRe} to realize $C\left( \mathbb{P}^{n}\left( \mathcal{T}%
\right) \right) $ as a concrete groupoid C*-algebra \cite{Rena}. Focusing on
the quantum complex line $\mathbb{P}^{1}\left( \mathcal{T}\right) $, we get
the C*-algebra structure of $C\left( \mathbb{P}^{1}\left( \mathcal{T}\right)
\right) $ analyzed and its $K$-groups computed. Furthermore, we get the
finitely generated left projective modules over $C\left( \mathbb{P}%
^{1}\left( \mathcal{T}\right) \right) $ classified up to isomorphism by
classifying the projections over $C\left( \mathbb{P}^{1}\left( \mathcal{T}%
\right) \right) $, i.e. in $M_{\infty}\left( C\left( \mathbb{P}^{1}\left( 
\mathcal{T}\right) \right) \right) $, up to unitary equivalence, and
explicitly identify the quantum line bundles $L_{k}$ among the classified
projections, showing that these modules $L_{k}$ do exhaust all rank-one
projections over $C\left( \mathbb{P}^{1}\left( \mathcal{T}\right) \right) $.
(After posting this work at ArXiv, a revised version of \cite{HaPaSiZi}
appeared at arXiv containing a computation of $K$-groups of $C\left( \mathbb{%
P}^{n}\left( \mathcal{T}\right) \right) $ for all $n$ \cite{HaNePaSiZi}.)

\section{Quantum projective spaces as groupoid C*-algebras}

Taking the groupoid approach to C*-algebras initiated by Renault \cite{Rena}
and popularized by the work of Curto, Muhly, and Renault \cite{CuMu,MuRe},
we give a description of the C*-algebras $C\left( \mathbb{S}%
_{H}^{2n-1}\right) $ and $C\left( \mathbb{P}^{n-1}\left( \mathcal{T}\right)
\right) $ of \cite{HaPaSiZi} as some concrete groupoid C*-algebras. We refer
to \cite{Rena,MuRe} for the concepts and theory of groupoid C*-algebras used
freely in the following discussion.

Let $\left. \left( \mathbb{Z}^{n}\ltimes\overline{\mathbb{Z}}^{n}\right)
\right\vert _{\overline{\mathbb{Z}_{\geq}}^{n}}$ with $n>1$ be the
transformation group groupoid $\mathbb{Z}^{n}\ltimes\overline{\mathbb{Z}}^{n}
$ restricted to the positive \textquotedblleft cone\textquotedblright\ $%
\overline{\mathbb{Z}_{\geq}}^{n}$ where $\overline{\mathbb{Z}}:=\mathbb{Z}%
\cup\left\{ \infty\right\} $ carries the standard topology and $\mathbb{Z}%
^{n}$ acts on $\overline{\mathbb{Z}}^{n}$ componentwise in the canonical
way. From the invariant open subset $\mathbb{Z}_{\geq}^{n}$ of the unit
space $\overline{\mathbb{Z}_{\geq}}^{n}$ of $\left. \left( \mathbb{Z}%
^{n}\ltimes\overline{\mathbb{Z}}^{n}\right) \right\vert _{\overline{\mathbb{Z%
}_{\geq}}^{n}}$, we get the short exact sequence of C*-algebras%
\begin{equation*}
0\rightarrow C^{\ast}\left( \left. \left( \mathbb{Z}^{n}\ltimes \overline{%
\mathbb{Z}}^{n}\right) \right\vert _{\mathbb{Z}_{\geq}^{n}}\right)
\rightarrow C^{\ast}\left( \left. \left( \mathbb{Z}^{n}\ltimes \overline{%
\mathbb{Z}}^{n}\right) \right\vert _{\overline{\mathbb{Z}_{\geq}}%
^{n}}\right) \rightarrow C^{\ast}\left( \left. \left( \mathbb{Z}^{n}\ltimes%
\overline{\mathbb{Z}}^{n}\right) \right\vert _{\overline {\mathbb{Z}_{\geq}}%
^{n}\backslash\mathbb{Z}_{\geq}^{n}}\right) \rightarrow0 
\end{equation*}
and furthermore since the open invariant set $\mathbb{Z}_{\geq}^{n}$ is
dense in the unit space, it induces a faithful representation $\pi$ of $%
C^{\ast }\left( \left. \left( \mathbb{Z}^{n}\ltimes\overline{\mathbb{Z}}%
^{n}\right) \right\vert _{\overline{\mathbb{Z}_{\geq}}^{n}}\right) $ on $%
\ell^{2}\left( \mathbb{Z}_{\geq}^{n}\right) $. From the groupoid isomorphism 
\begin{equation*}
\left. \left( \mathbb{Z}^{n}\ltimes\overline{\mathbb{Z}}^{n}\right)
\right\vert _{\overline{\mathbb{Z}_{\geq}}^{n}}\ \cong\ \times^{n}\left.
\left( \mathbb{Z}\ltimes\overline{\mathbb{Z}}\right) \right\vert _{\overline{%
\mathbb{Z}_{\geq}}}
\end{equation*}
and the C*-algebra isomorphism $C^{\ast}\left( \left. \left( \mathbb{Z}%
\ltimes\overline{\mathbb{Z}}\right) \right\vert _{\overline{\mathbb{Z}_{\geq
}}}\right) \cong\mathcal{T}$ for the Toeplitz C*-algebra $\mathcal{T}$ with $%
C^{\ast}\left( \left. \left( \mathbb{Z}\ltimes\overline{\mathbb{Z}}\right)
\right\vert _{\mathbb{Z}_{\geq}}\right) \cong\mathcal{K}\left(
\ell^{2}\left( \mathbb{Z}_{\geq}\right) \right) $, we get 
\begin{equation*}
C^{\ast}\left( \left. \left( \mathbb{Z}^{n}\ltimes\overline{\mathbb{Z}}%
^{n}\right) \right\vert _{\overline{\mathbb{Z}_{\geq}}^{n}}\right)
\cong\otimes^{n}\mathcal{T}
\end{equation*}
with $C^{\ast}\left( \left. \left( \mathbb{Z}^{n}\ltimes\overline {\mathbb{Z}%
}^{n}\right) \right\vert _{\mathbb{Z}_{\geq}^{n}}\right) \cong\otimes^{n}%
\mathcal{K}\left( \ell^{2}\left( \mathbb{Z}_{\geq}\right) \right) $.

Since $C\left( \mathbb{S}_{H}^{2n-1}\right) \cong\left( \otimes ^{n}\mathcal{%
T}\right) /\left( \otimes^{n}\mathcal{K}\right) $ by (A.2) of \cite{HaPaSiZi}%
, we have 
\begin{equation*}
C\left( \mathbb{S}_{H}^{2n-1}\right) \cong C^{\ast}\left( \left. \left( 
\mathbb{Z}^{n}\ltimes\overline{\mathbb{Z}}^{n}\right) \right\vert _{%
\overline{\mathbb{Z}_{\geq}}^{n}}\right) \diagup C^{\ast}\left( \left.
\left( \mathbb{Z}^{n}\ltimes\overline{\mathbb{Z}}^{n}\right) \right\vert _{%
\mathbb{Z}_{\geq}^{n}}\right) \cong C^{\ast}\left( \mathfrak{G}_{n}\right) 
\end{equation*}
realized as the groupoid C*-algebra of the concrete groupoid 
\begin{equation*}
\mathfrak{G}_{n}:=\left. \left( \mathbb{Z}^{n}\ltimes\overline{\mathbb{Z}}%
^{n}\right) \right\vert _{\overline{\mathbb{Z}_{\geq}}^{n}\backslash \mathbb{%
Z}_{\geq}^{n}}\ \ . 
\end{equation*}

Next we note that the $U\left( 1\right) $-action on $C\left( \mathbb{S}%
_{H}^{2n-1}\right) \cong C^{\ast}\left( \mathfrak{G}_{n}\right) $ considered
in \cite{HaPaSiZi} is implemented by the multiplication operator 
\begin{equation*}
U_{\zeta}:f\in C_{c}\left( \mathfrak{G}_{n}\right) \mapsto h_{\zeta}f\in
C_{c}\left( \mathfrak{G}_{n}\right) 
\end{equation*}
for $\zeta\in U\left( 1\right) \equiv\mathbb{T}$ where 
\begin{equation*}
h_{\zeta}:\left( m,p\right) \in\mathfrak{G}_{n}\subset\mathbb{Z}^{n}\ltimes%
\overline{\mathbb{Z}}^{n}\mapsto\zeta^{\Sigma m}\in\mathbb{T}\text{\ \ with }%
\Sigma m:=\sum_{i=1}^{n}m_{i}
\end{equation*}
is a groupoid character of $\mathfrak{G}_{n}$ and hence $U_{\zeta}$ is an
automorphism of $C^{\ast}\left( \mathfrak{G}_{n}\right) $, since $h_{\zeta
}h_{\zeta^{\prime}}=h_{\zeta\zeta^{\prime}}$ and $U_{\zeta}\left(
w_{i}\right) =\zeta w_{i}$ for the generators $w_{i}$ of $C\left( \mathbb{S}%
_{H}^{2n-1}\right) $ \cite{HaPaSiZi} identified with the characteristic
function $\chi_{A_{i}}\in C_{c}\left( \mathfrak{G}_{n}\right) $ of 
\begin{equation*}
A_{i}:=\left\{ \left( e_{i},p\right) :p\in\overline{\mathbb{Z}_{\geq}}%
^{n}\backslash\mathbb{Z}_{\geq}^{n}\right\} \subset\mathfrak{G}_{n}\subset%
\mathbb{Z}^{n}\ltimes\overline{\mathbb{Z}}^{n}. 
\end{equation*}

The C*-algebra $C\left( \mathbb{P}^{n-1}\left( \mathcal{T}\right) \right) $
of the quantum complex projective space studied in \cite{HaKaZi,HaPaSiZi} is
then isomorphic to the $U\left( 1\right) $-invariant C*-subalgebra $C^{\ast
}\left( \mathfrak{G}_{n}\right) ^{U\left( 1\right) }$ of $C^{\ast}\left( 
\mathfrak{G}_{n}\right) $, which can be realized as the groupoid C*-algebra $%
C^{\ast}\left( \left( \mathfrak{G}_{n}\right) _{0}\right) $ of the
subgroupoid $\left( \mathfrak{G}_{n}\right) _{0}$ of $\mathfrak{G}_{n}$,
where 
\begin{equation*}
\left( \mathfrak{G}_{n}\right) _{k}:=\left\{ \left( m,p\right) \in\mathfrak{G%
}_{n}:\Sigma m=k\right\} 
\end{equation*}
for $k\in\mathbb{Z}$. Furthermore, $C^{\ast}\left( \mathfrak{G}_{n}\right) $
becomes a graded algebra $\oplus_{k\in\mathbb{Z}}\overline{C_{c}\left(
\left( \mathfrak{G}_{n}\right) _{k}\right) }$ with the component $\overline{%
C_{c}\left( \left( \mathfrak{G}_{n}\right) _{k}\right) }$ being the quantum
line bundle $C\left( \mathbb{S}_{H}^{2n-1}\right) _{k}$ \cite%
{HaKaZi,HaPaSiZi} of degree $k$ over the quantum space $\mathbb{P}%
^{n-1}\left( \mathcal{T}\right) $.

As shown in \cite{HaPaSiZi}, the case of $\mathbb{P}^{n-1}\left( \mathcal{T}%
\right) $ with $n=2$ plays a crucially important role in the study of the
quantum line bundle $C\left( \mathbb{S}_{H}^{2n-1}\right) _{k}$ in general,
so we focus on the case of $\mathbb{P}^{1}\left( \mathcal{T}\right) $ in the
remaining part of this paper, while leaving the higher-dimensional cases to
a subsequent paper.

\section{$K$-groups of quantum projective line}

In the case of $n=2$, the groupoid $\mathcal{G}:=\mathfrak{G}%
_{2}\equiv\left. \left( \mathbb{Z}^{2}\ltimes\overline{\mathbb{Z}}%
^{2}\right) \right\vert _{\overline{\mathbb{Z}_{\geq}}^{2}\backslash\mathbb{Z%
}_{\geq}^{2}}$ has the unit space 
\begin{equation*}
\mathcal{G}^{\left( 0\right) }=\overline{\mathbb{Z}_{\geq}}^{2}\backslash%
\mathbb{Z}_{\geq}^{2}=\left( \overline{\mathbb{Z}_{\geq}}\times\left\{
\infty\right\} \right) \cup\left( \left\{ \infty\right\} \times\overline{%
\mathbb{Z}_{\geq}}\right) 
\end{equation*}
and consists of points $\left( m,l,p,q\right) $ with $\left( m,l\right) \in%
\mathbb{Z}^{2}$ and $\left( p,q\right) \in\mathcal{G}^{\left( 0\right) }$
such that $\left( m+p,l+q\right) \in\mathcal{G}^{\left( 0\right) }$ where $%
m+\infty=\infty$ for any $m\in\mathbb{Z}$ is understood.

The previous discussion realizes the C*-algebra $C\left( \mathbb{P}%
^{1}\left( \mathcal{T}\right) \right) $ of the quantum projective line $%
\mathbb{P}^{1}\left( \mathcal{T}\right) $ as a groupoid C*-algebra $%
C^{\ast}\left( \mathcal{G}_{0}\right) $ where the subgroupoid 
\begin{equation*}
\mathcal{G}_{0}:=\left\{ \left( n,-n,p,q\right) :n\in\mathbb{Z}\text{ such
that }\left( p,q\right) ,\left( p+n,q-n\right) \in\mathcal{G}^{\left(
0\right) }\right\} \subset\mathcal{G}
\end{equation*}
shares the same unit space $\left( \mathcal{G}_{0}\right) ^{\left( 0\right)
}=\mathcal{G}^{\left( 0\right) }$ with $\mathcal{G}$.

Note that the open dense invariant subset $U:=\left( \mathbb{Z}%
_{\geq}\times\left\{ \infty\right\} \right) \sqcup\left( \left\{
\infty\right\} \times\mathbb{Z}_{\geq}\right) $ of $\mathcal{G}^{\left(
0\right) }$ consists of two disjoint free orbits $\mathbb{Z}%
_{\geq}\times\left\{ \infty\right\} $ and $\left\{ \infty\right\} \times%
\mathbb{Z}_{\geq}$ of $\mathcal{G}_{0}$, from which we get a faithful
representation $\pi$ of $C^{\ast}\left( \mathcal{G}_{0}\right) \equiv
C\left( \mathbb{P}^{1}\left( \mathcal{T}\right) \right) $ on the Hilbert
space 
\begin{equation*}
\mathcal{H}:=\ell^{2}\left( U\right) \cong\ell^{2}\left( \mathbb{Z}_{\geq
}\right) \oplus\ell^{2}\left( \mathbb{Z}_{\geq}\right) 
\end{equation*}
such that $\pi\left( \delta_{\left( n,-n,p,q\right) }\right) $ for any $%
\left( n,-n,p,q\right) \in\mathcal{G}_{0}$ with $\left( p,q\right) \in U$ is
the partial isometry sending $\delta_{\left( p,q\right) }\in\ell ^{2}\left(
U\right) $ to $\delta_{\left( p+n,q-n\right) }\in\ell ^{2}\left( U\right) $
and all other $\delta_{\left( p^{\prime},q^{\prime }\right)
}\in\ell^{2}\left( U\right) $ to $0$.

The open subgroupoid 
\begin{equation*}
\mathcal{G}_{0}|_{U}=\left\{ \left( n,-n,p,\infty\right) ,\left(
n,-n,\infty,q\right) :n\in\mathbb{Z}\text{ such that }p,q,p+n,q-n\in \mathbb{%
Z}_{\geq}\right\} 
\end{equation*}
of $\mathcal{G}_{0}$ is isomorphic to the disjoint union $\mathfrak{K}%
_{+}\sqcup\mathfrak{K}_{-}$ of two copies of the groupoid $\mathfrak{K}%
:=\left. \left( \mathbb{Z}\ltimes\mathbb{Z}\right) \right\vert _{\mathbb{Z}%
_{\geq}}$ under the map $\left( n,-n,p,\infty\right) \mapsto\left(
n,p\right) \in\mathfrak{K}_{+}$ and $\left( n,-n,\infty ,q\right)
\mapsto\left( -n,q\right) \in\mathfrak{K}_{-}$. Thus 
\begin{equation*}
C^{\ast}\left( \mathcal{G}_{0}|_{U}\right) \cong C^{\ast}\left( \mathfrak{K}%
_{+}\right) \oplus C^{\ast}\left( \mathfrak{K}_{-}\right) \cong\mathcal{K}%
\left( \ell^{2}\left( \mathbb{Z}_{\geq}\right) \right) \oplus\mathcal{K}%
\left( \ell^{2}\left( \mathbb{Z}_{\geq}\right) \right) . 
\end{equation*}
With $\mathcal{G}^{\left( 0\right) }\backslash U=\left\{ \left(
\infty,\infty\right) \right\} $ and $\mathcal{G}_{0}|_{\left\{ \left(
\infty,\infty\right) \right\} }=\left\{ \left( n,-n,\infty,\infty\right)
:n\in\mathbb{Z}\right\} $ isomorphic to the group $\mathbb{Z}$, we have the
short exact sequence

\begin{equation*}
0\rightarrow C^{\ast}\left( \left. \mathcal{G}_{0}\right\vert _{U}\right)
\cong\mathcal{K}\oplus\mathcal{K}\overset{\iota}{\rightarrow}C^{\ast}\left( 
\mathcal{G}_{0}\right) \overset{\sigma}{\rightarrow}C^{\ast}\left( \left. 
\mathcal{G}_{0}\right\vert _{\left\{ \left( \infty,\infty\right) \right\}
}\right) \cong C\left( \mathbb{T}\right) \rightarrow0, 
\end{equation*}
where $C^{\ast}\left( \left. \mathcal{G}_{0}\right\vert _{U}\right) \cong%
\mathcal{K}\oplus\mathcal{K}$ under the representation $\pi$ with $\mathcal{K%
}\equiv\mathcal{K}\left( \ell^{2}\left( \mathbb{Z}_{\geq}\right) \right) $,
and $\delta_{\left( 1,-1,\infty,\infty\right) }\in C^{\ast }\left( \left. 
\mathcal{G}_{0}\right\vert _{\left\{ \left( \infty ,\infty\right) \right\}
}\right) $ is identified with 
\begin{equation*}
z:=\func{id}_{\mathbb{T}}\in C\left( \mathbb{T}\right) . 
\end{equation*}

In the induced 6-term exact sequence%
\begin{equation*}
\begin{array}{rcccccl}
\mathbb{Z}\left[ e_{11}\right] \oplus\mathbb{Z}\left[ e_{11}\right] \cong & 
K_{0}\left( \mathcal{K}\oplus\mathcal{K}\right) & \overset{K_{0}\left(
\iota\right) }{\rightarrow} & K_{0}\left( C^{\ast}\left( \mathcal{G}%
_{0}\right) \right) & \overset{K_{0}\left( \sigma\right) }{\rightarrow} & 
K_{0}\left( C\left( \mathbb{T}\right) \right) & \cong\mathbb{Z} \\ 
& \uparrow\eta &  &  &  & \downarrow\varepsilon &  \\ 
\mathbb{Z}\left[ z\right] \cong & K_{1}\left( C\left( \mathbb{T}\right)
\right) & \overset{K_{1}\left( \sigma\right) }{\leftarrow} & K_{1}\left(
C^{\ast}\left( \mathcal{G}_{0}\right) \right) & \overset{K_{1}\left(
\iota\right) }{\leftarrow} & K_{1}\left( \mathcal{K}\oplus\mathcal{K}\right)
& =0,%
\end{array}
\end{equation*}
the homomorphism $K_{0}\left( \sigma\right) $ is clearly surjective, and we
claim that the index homomorphism $\eta$ sends $\left[ z\right] $ to $\left(
-\left[ e_{11}\right] \right) \oplus\left[ e_{11}\right] $, where $e_{11}$
is the standard matrix unit and $z\equiv\func{id}_{\mathbb{T}}\in
GL_{1}\left( C\left( \mathbb{T}\right) \right) $. Indeed $z$ lifts via $%
\sigma$ to the characteristic function $\chi_{W}\in C^{\ast }\left( \mathcal{%
G}_{0}\right) $ of the set 
\begin{equation*}
W:=\left\{ \left( 1,-1,p,\infty\right) :p\geq0\right\} \cup\left\{ \left(
1,-1,\infty,q\right) :q\geq1\right\} , 
\end{equation*}
and $\pi\left( \chi_{W}\right) =\mathcal{S}\oplus\mathcal{S}^{\ast}$ a
partial isometry with kernel projection $0\oplus e_{11}$ and cokernel
projection $e_{11}\oplus0$, where $\mathcal{S}$ is the (forward) unilateral
shift. Hence 
\begin{equation*}
\eta\left( \left[ z\right] \right) =\left[ 0\oplus e_{11}\right] -\left[
e_{11}\oplus0\right] \in K_{0}\left( \mathcal{K}\oplus \mathcal{K}\right) . 
\end{equation*}
(It is understood that the index homomorphism $\eta$ used here may be
different by a $\pm$-sign from the one used by other authors.)

Thus we get $K_{0}\left( \iota\right) \left( \left[ 0\oplus e_{11}\right] -%
\left[ e_{11}\oplus0\right] \right) =0$ in $K_{0}\left( C^{\ast}\left( 
\mathcal{G}_{0}\right) \right) $, and hence $\left[ e_{11}\oplus0\right] =%
\left[ 0\oplus e_{11}\right] $ in $K_{0}\left( C^{\ast}\left( \mathcal{G}%
_{0}\right) \right) $. A simple diagram chase concludes that $K_{0}\left(
C^{\ast}\left( \mathcal{G}_{0}\right) \right) \cong\mathbb{Z}\left[
e_{11}\oplus0\right] \oplus\mathbb{Z}\left[ \tilde {I}\right] $ for the
identity element $\tilde{I}$ of $C^{\ast}\left( \mathcal{G}_{0}\right) $,
while $K_{1}\left( C^{\ast}\left( \mathcal{G}_{0}\right) \right) =0$.
Furthermore 
\begin{align*}
K_{0}\left( \iota\right) & :m\left[ e_{11}\right] \oplus l\left[ e_{11}%
\right] \in K_{0}\left( \mathcal{K}\oplus\mathcal{K}\right) \cong\mathbb{Z}%
\oplus\mathbb{Z} \\
& \mapsto\left( m+l\right) \left[ e_{11}\oplus0\right] \oplus0\left[ \tilde{I%
}\right] \in K_{0}\left( C^{\ast}\left( \mathcal{G}_{0}\right) \right) \cong%
\mathbb{Z}\oplus\mathbb{Z}.
\end{align*}

We summarize as follows.

\textbf{Theorem 1}. \textit{For the quantum complex projective space }$%
\mathbb{P}^{1}\left( \mathcal{T}\right) $\textit{, there is a short exact
sequence of C*-algebras decomposing its algebra }$C\left( \mathbb{P}%
^{1}\left( \mathcal{T}\right) \right) $\textit{\ as }%
\begin{equation*}
0\rightarrow\mathcal{K}\oplus\mathcal{K}\overset{\iota}{\rightarrow}C\left( 
\mathbb{P}^{1}\left( \mathcal{T}\right) \right) \overset{\sigma }{\rightarrow%
}C\left( \mathbb{T}\right) \rightarrow0, 
\end{equation*}
\textit{and its }$K$\textit{-groups coincide with those of its classical
counterpart, i.e. }%
\begin{equation*}
K_{0}\left( C\left( \mathbb{P}^{1}\left( \mathcal{T}\right) \right) \right)
\cong\mathbb{Z}\oplus\mathbb{Z}\text{\ \ \textit{and}\ \ }K_{1}\left(
C\left( \mathbb{P}^{1}\left( \mathcal{T}\right) \right) \right) =0. 
\end{equation*}

\section{Classification of projections over quantum projective line}

In the following, we denote by $M_{\infty}\left( \mathcal{A}\right) $ the
direct limit (or the union as sets) of the increasing sequence of matrix
algebras $M_{n}\left( \mathcal{A}\right) $ over $\mathcal{A}$ with the
canonical inclusion $M_{n}\left( \mathcal{A}\right) \subset M_{n+1}\left( 
\mathcal{A}\right) $ identifying $x\in M_{n}\left( \mathcal{A}\right) $ with 
$x\boxplus0\in M_{n+1}\left( \mathcal{A}\right) $ for any algebra $\mathcal{A%
}$, where $\boxplus$ denotes the standard diagonal concatenation of two
matrices. So the size of an element in $M_{\infty}\left( \mathcal{A}\right) $
can be taken arbitrarily large. We also use $U_{\infty}\left( \mathcal{A}%
\right) $ to denote the direct limit of the unitary groups $U_{n}\left( 
\mathcal{A}\right) \subset M_{n}\left( \mathcal{A}\right) $ for a unital
C*-algebra $\mathcal{A}$ with $U_{n}\left( \mathcal{A}\right) $ embedded in $%
U_{n+1}\left( \mathcal{A}\right) $ by identifying $x\in U_{n}\left( \mathcal{%
A}\right) $ with $x\boxplus1\in U_{n+1}\left( \mathcal{A}\right) $.

Before proceeding with the classification problem, we briefly recall the
relation between projections over a C*-algebra $\mathcal{A}$ and finitely
generated left projective modules over $\mathcal{A}$, and between them and $K
$-theory.

By a projection over a unital C*-algebra $\mathcal{A}$, we mean a
self-adjoint idempotent in $M_{\infty}\left( \mathcal{A}\right) $. Two
projections $P,Q\in M_{n}\left( \mathcal{A}\right) $ are called unitarily
equivalent if there exists a unitary $U\in M_{N}\left( \mathcal{A}\right) $
with $N\geq n$ such that $UPU^{-1}=Q$. Each projection $P\in M_{n}\left( 
\mathcal{A}\right) $ over $\mathcal{A}$ defines a finitely generated left
projective module $\mathcal{A}^{n}P$ over $\mathcal{A}$ where elements of $%
\mathcal{A}^{n}$ are viewed as row vectors. The mapping $P\mapsto\mathcal{A}%
^{n}P$ induces a bijective correspondence between the unitary equivalence
classes of projections over $\mathcal{A}$ and the isomorphism classes of
finitely generated left projective modules over $\mathcal{A}$ \cite{Blac}.

Two finitely generated projective left modules $E,F$ over $\mathcal{A}$ are
called stably isomorphic if they become isomorphic after being augmented by
the same finitely generated free $\mathcal{A}$-module, i.e. $E\oplus 
\mathcal{A}^{k}\cong F\oplus\mathcal{A}^{k}$ for some $k\geq0$.
Correspondingly, two projections $P$ and $Q$ are called stably equivalent if 
$P\boxplus I_{k}$ and $Q\boxplus I_{k}$ are unitarily equivalent for some
identity matrix $I_{k}$. The $K_{0}$-group $K_{0}\left( \mathcal{A}\right) $
classifies projections over $\mathcal{A}$ up to stable equivalence. The
classification of projections over a C*-algebra up to unitary equivalence,
appearing as the cancellation problem, was popularized by Rieffel's
pioneering work \cite{Ri:dsr,Ri:ct} and is in general an interesting but
difficult question.

The set of all unitary equivalence classes of projections over a C*-algebra $%
\mathcal{A}$ is an abelian monoid $\mathfrak{P}\left( \mathcal{A}\right) $
with its binary operation provided by the diagonal sum $\boxplus$ of
projections. The image of the canonical homomorphism from $\mathfrak{P}%
\left( \mathcal{A}\right) $ into $K_{0}\left( \mathcal{A}\right) $ is the
so-called positive cone of $K_{0}\left( \mathcal{A}\right) $.

In the following, we use $\tilde{I}$ to denote the multiplicative unit of
the unital C*-algebra $\left( \mathcal{K}\oplus\mathcal{K}\right)
^{+}\subset C^{\ast}\left( \mathcal{G}_{0}\right) $ where $\mathcal{A}^{+}$
denotes the unitization of $\mathcal{A}$, and $\tilde{I}_{n}$ to denote the
identity matrix in $M_{n}\left( \left( \mathcal{K}\oplus\mathcal{K}\right)
^{+}\right) $, while 
\begin{equation*}
P_{m}:=\sum_{i=1}^{m}e_{ii}\in M_{m}\left( \mathbb{C}\right) \subset 
\mathcal{K}
\end{equation*}
denotes the standard $m\times m$ identity matrix in $M_{m}\left( \mathbb{C}%
\right) \subset\mathcal{K}$ for any integer $m\geq0$ (with $M_{0}\left( 
\mathbb{C}\right) \equiv\left\{ 0\right\} $ and $P_{0}\equiv0$ understood).
We also use the notation 
\begin{equation*}
P_{-m}:=I-P_{m}\in\mathcal{K}^{+}
\end{equation*}
for integers $m>0$, where $I$ is the identity operator canonically contained
in $\mathcal{K}^{+}$, and symbolically adopt the notation 
\begin{equation*}
P_{-0}\equiv I-P_{0}=I\neq P_{0}. 
\end{equation*}
Furthermore by abuse of notation, we take 
\begin{equation*}
P_{-m}\oplus P_{-l}:=\tilde{I}-\left( P_{m}\oplus P_{l}\right) \in\left( 
\mathcal{K}\oplus\mathcal{K}\right) ^{+}
\end{equation*}
if $m,l\geq0$. Note that $P_{m}\oplus P_{l}\notin\left( \mathcal{K}\oplus%
\mathcal{K}\right) ^{+}$ if $m$ and $l$ are of strictly opposite $\pm$-sings.

Let $\alpha\in M_{\infty}\left( C^{\ast}\left( \mathcal{G}_{0}\right)
\right) $ be a projection. Since projections in $M_{\infty}\left( C\left( 
\mathbb{T}\right) \right) $ are classified up to unitary equivalence as the
constant functions on $\mathbb{T}$ with an identity matrix $I_{n}\in
M_{n}\left( \mathbb{C}\right) $ as the value for some $n\in\mathbb{Z}_{\geq }
$ (and hence $K_{0}\left( C\left( \mathbb{T}\right) \right) =\mathbb{Z}$), $%
\alpha$ is unitarily equivalent over $C^{\ast}\left( \mathcal{G}_{0}\right) $
to some projection $\beta\in M_{N}\left( C^{\ast}\left( \mathcal{G}%
_{0}\right) \right) $ with $\sigma\left( \beta\right) =I_{n}$ for some $%
n\geq0$ and a suitably large size $N\geq n$. It is easy to see that $n$
depends only on $\alpha$, and we call $n$ the rank of $\alpha$.

So in the following, we concentrate on classifying projections $\alpha$ over 
$C^{\ast}\left( \mathcal{G}_{0}\right) $ with $\sigma\left( \alpha\right)
=I_{n}$ and $\alpha\in M_{N}\left( C^{\ast}\left( \mathcal{G}_{0}\right)
\right) $ for some $N\geq n$.

Now since $\sigma\left( \alpha-\tilde{I}_{n}\right) =I_{n}-I_{n}=0$, 
\begin{equation*}
\alpha-\tilde{I}_{n}\in M_{N}\left( \mathcal{K}\oplus\mathcal{K}\right)
\equiv M_{N}\left( \mathcal{K}\right) \oplus M_{N}\left( \mathcal{K}\right) 
\end{equation*}
which can be approximated by elements in $M_{N}\left( M_{k}\left( \mathbb{C}%
\right) \right) \oplus M_{N}\left( M_{k}\left( \mathbb{C}\right) \right) $.
So we can replace $\alpha$ by a unitarily equivalent projection $\tilde{I}%
_{n}+x$ for some $x$ in 
\begin{equation*}
M_{N}\left( M_{k}\left( \mathbb{C}\right) \right) \oplus M_{N}\left(
M_{k}\left( \mathbb{C}\right) \right) \subset M_{N}\left( \pi\left(
C_{c}\left( \left. \mathcal{G}_{0}\right\vert _{U}\right) \right) \right)
\subset M_{n}\left( \left( \mathcal{K}\oplus\mathcal{K}\right) ^{+}\right) 
\end{equation*}
with a suitably large $k$. Let $I_{n}^{\prime}$ be the identity element of 
\begin{equation*}
M_{n}\left( M_{k}\left( \mathbb{C}\right) \right) \oplus M_{n}\left(
M_{k}\left( \mathbb{C}\right) \right) \subset M_{n}\left( \left( \mathcal{K}%
\oplus\mathcal{K}\right) ^{+}\right) . 
\end{equation*}
Then since $I_{n}^{\prime}+x\in M_{N}\left( M_{k}\left( \mathbb{C}\right)
\right) \oplus M_{N}\left( M_{k}\left( \mathbb{C}\right) \right) $ is
unitarily equivalent over $M_{k}\left( \mathbb{C}\right) \oplus M_{k}\left( 
\mathbb{C}\right) $ to $P_{m}\oplus P_{l}$ for some $0\leq m,l\leq Nk$ where 
$P_{m}\oplus P_{l}$ is the identity element of 
\begin{equation*}
M_{m}\left( \mathbb{C}\right) \oplus M_{l}\left( \mathbb{C}\right) \subset
M_{Nk}\left( \mathbb{C}\right) \oplus M_{Nk}\left( \mathbb{C}\right) \equiv
M_{N}\left( M_{k}\left( \mathbb{C}\right) \right) \oplus M_{N}\left(
M_{k}\left( \mathbb{C}\right) \right) , 
\end{equation*}
we have $\tilde{I}_{n}+x$ unitarily equivalent over $\left( M_{k}\left( 
\mathbb{C}\right) \oplus M_{k}\left( \mathbb{C}\right) \right)
^{+}\subset\pi\left( C_{c}\left( \left. \mathcal{G}_{0}\right\vert
_{U}\right) \right) ^{+}$ to 
\begin{equation*}
\left( \tilde{I}_{n}-I_{n}^{\prime}\right) +\left( P_{m}\oplus P_{l}\right)
\in M_{N}\left( \left( M_{k}\left( \mathbb{C}\right) \oplus M_{k}\left( 
\mathbb{C}\right) \right) ^{+}\right) \subset M_{N}\left( \left( \mathcal{K}%
\oplus\mathcal{K}\right) ^{+}\right) 
\end{equation*}
by the canonical embedding of $U_{\infty}\left( \mathcal{A}\right) $ in $%
U_{\infty}\left( \mathcal{A}^{+}\right) $ for any unital C*-algebra $%
\mathcal{A}$. Note that $\left( \tilde{I}_{n}-I_{n}^{\prime}\right) +\left(
P_{m}\oplus P_{l}\right) $ can be expressed in the form 
\begin{equation*}
\text{(*)\ \ }\left( P_{m_{1}}\oplus P_{l_{1}}\right) \boxplus\cdots
\boxplus\left( P_{m_{N}}\oplus P_{l_{N}}\right) \in M_{N}\left( \left(
M_{k}\left( \mathbb{C}\right) \oplus M_{k}\left( \mathbb{C}\right) \right)
^{+}\right) 
\end{equation*}
for some $m_{i},l_{i}\in\mathbb{Z}$ with $\left\vert m_{i}\right\vert
,\left\vert l_{i}\right\vert \leq k$, and since $\sigma\left( \alpha\right)
=\sigma\left( \tilde{I}_{n}\right) =I_{n}$, we have $m_{i},l_{i}\leq0$
(viewing $P_{m_{i}}\oplus P_{l_{1}}$ as $\tilde{I}-\left( P_{\left\vert
m_{i}\right\vert }\oplus P_{\left\vert l_{1}\right\vert }\right) $) for $%
i\leq n$ and $m_{i},l_{i}\geq0$ for $i>n$.

It remains to classify projections $\alpha\in M_{N}\left( C^{\ast}\left( 
\mathcal{G}_{0}\right) \right) $ of the form (*) up to unitary equivalence
over $C^{\ast}\left( \mathcal{G}_{0}\right) $.

When the rank $n$ is $0$, we have $m_{i},l_{i}\geq0$ for all $i$ in (*), and 
$P_{m_{i}},P_{l_{i}}\in M_{k}\left( \mathbb{C}\right) $. With $%
P_{m_{i}},P_{l_{i}}$ viewed as elements in $M_{Nk}\left( \mathbb{C}\right)
\supset M_{k}\left( \mathbb{C}\right) $, the projections $P_{m_{1}}\boxplus
\cdots\boxplus P_{m_{N}}$ and $P_{l_{1}}\boxplus\cdots\boxplus P_{l_{N}}$
lying in $M_{N}\left( M_{Nk}\left( \mathbb{C}\right) \right) $ with ranks
bounded by $Nk$ are unitarily equivalent over $M_{Nk}\left( \mathbb{C}%
\right) $ to $P_{m}\boxplus0\boxplus\cdots\boxplus0$ and $%
P_{l}\boxplus0\boxplus\cdots\boxplus0$ in $M_{N}\left( M_{Nk}\left( \mathbb{C%
}\right) \right) $ respectively, where $m:=\sum_{i}m_{i}\leq Nk$ and $%
l:=\sum_{i}l_{i}\leq Nk$. Hence $\left( P_{m_{1}}\oplus P_{l_{1}}\right)
\boxplus\cdots\boxplus\left( P_{m_{N}}\oplus P_{l_{N}}\right) $ is unitarily
equivalent over $\pi\left( C_{c}\left( \left. \mathcal{G}_{0}\right\vert
_{U}\right) \right) ^{+}$ to $P_{m}\oplus P_{l}\in M_{1}\left( \left( 
\mathcal{K}\oplus\mathcal{K}\right) ^{+}\right) \equiv\left( \mathcal{K}%
\oplus\mathcal{K}\right) ^{+}$.

On the other hand, if such projections $P_{m}\oplus P_{l}$ and $%
P_{m^{\prime}}\oplus P_{l^{\prime}}$ with $m,l,m^{\prime},l^{\prime}\in%
\mathbb{Z}_{\geq}$ are unitarily equivalent over $C^{\ast}\left( \mathcal{G}%
_{0}\right) \subset\mathcal{B}\left( \ell^{2}\left( \mathbb{Z}_{\geq}\right)
\right) \oplus\mathcal{B}\left( \ell^{2}\left( \mathbb{Z}_{\geq}\right)
\right) $, then their ranks must coincide, i.e. $m=m^{\prime}$ and $%
l=l^{\prime}$. Thus for the case of $n=0$, we get unitary equivalence
classes of projections over $C^{\ast}\left( \mathcal{G}_{0}\right) $
classified by $\left( m,l\right) \in\mathbb{Z}_{\geq}\times\mathbb{Z}_{\geq}$
as $P_{m}\oplus P_{l}$.

When the rank $n$ is strictly positive, we claim that the projections $%
\tilde{I}_{n}$, $\tilde{I}_{n}\boxplus\left( P_{m}\oplus P_{0}\right) $, and 
\begin{equation*}
\tilde{I}_{n-1}\boxplus\left( P_{-m}\oplus P_{-0}\right) \equiv\tilde {I}%
_{n-1}\boxplus\left( \tilde{I}-\left( P_{m}\oplus0\right) \right) 
\end{equation*}
with $m\in\mathbb{N}\equiv\mathbb{Z}_{>}$, give a complete list of unitary
equivalence classes of projections $\alpha$ with $\sigma\left( \alpha\right)
=I_{n}$.

First we observe that for $m,l,n\geq0$. 
\begin{equation*}
\left[ \tilde{I}_{n}\boxplus\left( P_{m}\oplus P_{l}\right) \right] =\left[ 
\tilde{I}_{n}\right] +K_{0}\left( \iota\right) \left( m\left[ e_{11}\right]
\oplus l\left[ e_{11}\right] \right) =\left( m+l\right) \left[ e_{11}\oplus0%
\right] \oplus n\left[ \tilde{I}\right] 
\end{equation*}
and 
\begin{equation*}
\left[ \tilde{I}_{n}\boxplus\left( P_{-m}\oplus P_{-l}\right) \right] =\left[
\tilde{I}_{n}\boxplus\left( \tilde{I}-\left( P_{m}\oplus P_{l}\right)
\right) \right] =\left[ \tilde{I}_{n}\right] +\left[ \tilde{I}-\left(
P_{m}\oplus P_{l}\right) \right] 
\end{equation*}%
\begin{equation*}
=\left[ \tilde{I}_{n}\right] +\left[ \tilde{I}\right] -\left[ P_{m}\oplus
P_{l}\right] =\left( n+1\right) \left[ \tilde{I}\right] -K_{0}\left(
\iota\right) \left( m\left[ e_{11}\right] \oplus l\left[ e_{11}\right]
\right) 
\end{equation*}%
\begin{equation*}
=-\left( m+l\right) \left[ e_{11}\oplus0\right] \oplus\left( n+1\right) %
\left[ \tilde{I}\right] 
\end{equation*}
in $K_{0}\left( C\left( \mathbb{P}^{1}\left( \mathcal{T}\right) \right)
\right) $. So the stable equivalence class over $C\left( \mathbb{P}%
^{1}\left( \mathcal{T}\right) \right) $ of a projection of the form $\tilde{I%
}_{n}\boxplus\left( P_{m}\oplus P_{l}\right) $ with $n,ml\geq0$ (so $m,l$
are integers not of opposite $\pm$-signs) is determined exactly by $m+l\in%
\mathbb{Z}$ and its rank ($n$ or $n+1$). In particular, for $n>0$, the
projections $\tilde{I}_{n}$, $\tilde{I}_{n}\boxplus\left( P_{m}\oplus
P_{0}\right) $, and $\tilde{I}_{n-1}\boxplus\left( P_{-m}\oplus
P_{-0}\right) $ with $m>0$ are mutually stably and hence unitarily
inequivalent.

It remains to show that any projection $\alpha$ of the form (*) with $n>0$
is unitarily equivalent to one of $\tilde{I}_{n}$, $\tilde{I}%
_{n}\boxplus\left( P_{m}\oplus P_{0}\right) $, and $\tilde{I}%
_{n-1}\boxplus\left( P_{-m}\oplus P_{-0}\right) $ with $m\in\mathbb{N}$.

Recall that $\sigma\left( \alpha\right) =I_{n}$ implies that $%
m_{i},l_{i}\leq0$ for $i\leq n$ and $m_{i},l_{i}\geq0$ for $i>n$. Since $%
\mathcal{K}\oplus\mathcal{K}\subset C\left( \mathbb{P}^{1}\left( \mathcal{T}%
\right) \right) $, using some (unitary) finite permutation matrices, we can
convert $\alpha$ to a unitarily equivalent projection $\beta$ of the form%
\begin{equation*}
\beta=\left( \boxplus^{n-1}\tilde{I}\right) \boxplus\left( P_{m^{\prime
\prime}}\oplus P_{l^{\prime\prime}}\right) \boxplus\left(
P_{m^{\prime}}\oplus P_{l^{\prime}}\right) \boxplus\left(
\boxplus^{N-n-1}0\right) \in M_{N}\left( \left( \mathcal{K}\oplus\mathcal{K}%
\right) ^{+}\right) 
\end{equation*}
with $m^{\prime\prime}=\sum_{i=1}^{n}m_{i}\leq0$, $l^{\prime\prime}=\sum
_{i=1}^{n}l_{i}\leq0$, $m^{\prime}=\sum_{i=n+1}^{N}m_{i}\geq0$, $l^{\prime
}=\sum_{i=n+1}^{N}l_{i}\geq0$, or for short 
\begin{equation*}
\beta=\tilde{I}_{n-1}\boxplus\left( P_{m^{\prime\prime}}\oplus P_{l^{\prime
\prime}}\right) \boxplus\left( P_{m^{\prime}}\oplus P_{l^{\prime}}\right)
\in M_{n+1}\left( \left( \mathcal{K}\oplus\mathcal{K}\right) ^{+}\right) , 
\end{equation*}
by swapping the largest (finite) identity diagonal blocks in $P_{m_{i}}$ and 
$P_{l_{i}}$ for $i>n+1$ with suitable disjoint diagonal zero blocks of $%
P_{m_{n+1}}$ and $P_{l_{n+1}}$ respectively, and by swapping the largest
(finite) diagonal zero blocks in $P_{m_{i}}$ and $P_{l_{i}}$ for $i<n$ with
suitable disjoint diagonal identity blocks of $P_{m_{n}}$ and $P_{l_{n}}$
respectively. Here it is understood that $m^{\prime\prime}$ and $l^{\prime
\prime}$ carry a negative sign and hence $P_{m^{\prime\prime}}\oplus
P_{l^{\prime\prime}}=\tilde{I}-\left( P_{\left\vert m^{\prime\prime
}\right\vert }\oplus P_{\left\vert l^{\prime\prime}\right\vert }\right) $.

By swapping a suitable (finite) diagonal zero block of $P_{l^{\prime\prime}}$
with a suitable (finite) identity block of $P_{l^{\prime}}$, we get $\beta$
unitarily equivalent to either 
\begin{equation*}
\tilde{I}_{n-1}\boxplus\left( P_{m^{\prime\prime}}\oplus P_{l^{\prime\prime
}+l^{\prime}}\right) \boxplus\left( P_{m^{\prime}}\oplus0\right) 
\end{equation*}
if $l^{\prime\prime}+l^{\prime}<0$, or to 
\begin{equation*}
\tilde{I}_{n-1}\boxplus\left( P_{m^{\prime\prime}}\oplus P_{-0}\right)
\boxplus\left( P_{m^{\prime}}\oplus P_{l^{\prime}+l^{\prime\prime}}\right) 
\end{equation*}
if $l^{\prime\prime}+l^{\prime}\geq0$.

With $\pi\left( \chi_{W}\right) =\mathcal{S}\oplus\mathcal{S}^{\ast}$ as
discussed earlier, conjugating $\tilde{I}_{n-1}\boxplus\left( P_{m^{\prime
\prime}}\oplus P_{l^{\prime\prime}+l^{\prime}}\right) \boxplus\left(
P_{m^{\prime}}\oplus0\right) $ or $\tilde{I}_{n-1}\boxplus\left(
P_{m^{\prime\prime}}\oplus P_{-0}\right) \boxplus\left( P_{m^{\prime}}\oplus
P_{l^{\prime}+l^{\prime\prime}}\right) $ by the unitary 
\begin{equation*}
\tilde{I}_{n-1}\boxplus\left( 
\begin{array}{cc}
\pi\left( \chi_{W}\right) ^{\left\vert l^{\prime\prime}+l^{\prime
}\right\vert } & \tilde{I}-\pi\left( \chi_{W}\right) ^{\left\vert
l^{\prime\prime}+l^{\prime}\right\vert }\left( \pi\left( \chi_{W}\right)
^{\ast}\right) ^{\left\vert l^{\prime\prime}+l^{\prime}\right\vert } \\ 
\tilde{I}-\left( \pi\left( \chi_{W}\right) ^{\ast}\right) ^{\left\vert
l^{\prime\prime}+l^{\prime}\right\vert }\pi\left( \chi_{W}\right)
^{\left\vert l^{\prime\prime}+l^{\prime}\right\vert } & \left( \pi\left(
\chi_{W}\right) ^{\ast}\right) ^{\left\vert l^{\prime\prime}+l^{\prime
}\right\vert }%
\end{array}
\right) \in M_{n+1}\left( C\left( \mathbb{P}^{1}\left( \mathcal{T}\right)
\right) \right) 
\end{equation*}
or its adjoint respectively converts each to the form 
\begin{equation*}
\gamma:=\tilde{I}_{n-1}\boxplus\left( P_{-j}\oplus P_{-0}\right)
\boxplus\left( P_{k}\oplus0\right) 
\end{equation*}
for some $j,k\geq0$ (up to swap of finite diagonal blocks in the first $%
\oplus$-summand).

Finally by swapping a suitable (finite) diagonal zero block of $P_{-j}$ with
a suitable (finite) identity block of $P_{k}$, we get $\gamma$ unitarily
equivalent to either $\tilde{I}_{n-1}\boxplus\left( P_{k-j}\oplus
P_{-0}\right) \boxplus0$ if $k-j<0$, or $\tilde{I}_{n}\boxplus\left(
P_{k-j}\oplus0\right) $ if $k-j\geq0$. Thus $\alpha$ is unitarily equivalent
over $C\left( \mathbb{P}^{1}\left( \mathcal{T}\right) \right) $ to either $%
\tilde{I}_{n-1}\boxplus\left( \tilde{I}-\left( P_{j-k}\oplus0\right) \right) 
$ if $k-j<0$, or $\tilde{I}_{n}\boxplus\left( P_{k-j}\oplus0\right) $ if $%
k-j\geq0$ as wanted.

Now we summarize what we have found.

\textbf{Theorem 2}. \textit{The abelian monoid }$\mathfrak{P}\left( C\left( 
\mathbb{P}^{1}\left( \mathcal{T}\right) \right) \right) $\textit{\ of
unitary equivalence classes of projections over }$C\left( \mathbb{P}%
^{1}\left( \mathcal{T}\right) \right) $\textit{\ consists of (the
representatives) }$P_{m}\oplus P_{l}$\textit{, }$\tilde{I}_{n}\boxplus\left(
P_{j}\oplus0\right) $\textit{, and }$\tilde{I}_{n-1}\boxplus\left( \tilde {I}%
-\left( P_{k}\oplus0\right) \right) $\textit{\ for }$m,l,j\in \mathbb{Z}%
_{\geq}$\textit{\ and }$n,k\in\mathbb{Z}_{>}$\textit{, where }$\tilde{I}$%
\textit{\ is the identity of }$\left( \mathcal{K}\oplus \mathcal{K}\right)
^{+}\subset C\left( \mathbb{P}^{1}\left( \mathcal{T}\right) \right) $\textit{%
\ and }$P_{k}$\textit{\ is the identity element of }$M_{k}\left( \mathbb{C}%
\right) \subset\mathcal{K}$\textit{, with its binary operation }$\cdot$%
\textit{\ specified by}%
\begin{equation*}
\left\{ 
\begin{array}{lll}
\left( P_{m}\oplus P_{l}\right) \cdot\left( \tilde{I}_{n}\boxplus\left(
P_{j}\oplus0\right) \right) =\tilde{I}_{n}\boxplus\left(
P_{m+l+j}\oplus0\right) , &  &  \\ 
\left( \tilde{I}_{n}\boxplus\left( P_{j}\oplus0\right) \right) \cdot\left( 
\tilde{I}_{n^{\prime}-1}\boxplus\left( \tilde{I}-\left( P_{k}\oplus0\right)
\right) \right) =\tilde{I}_{n+n^{\prime}}\boxplus\left( P_{j-k}\oplus0\right)
& \text{\textit{if }} & j\geq k, \\ 
\left( \tilde{I}_{n}\boxplus\left( P_{j}\oplus0\right) \right) \cdot\left( 
\tilde{I}_{n^{\prime}-1}\boxplus\left( \tilde{I}-\left( P_{k}\oplus0\right)
\right) \right) =\tilde{I}_{n+n^{\prime}-1}\boxplus\left( \tilde{I}-\left(
P_{k-j}\oplus0\right) \right) & \text{\textit{if }} & j<k, \\ 
\left( P_{m}\oplus P_{l}\right) \cdot\left( \tilde{I}_{n-1}\boxplus\left( 
\tilde{I}-\left( P_{k}\oplus0\right) \right) \right) =\tilde{I}%
_{n}\boxplus\left( P_{m+l-k}\oplus0\right) & \text{\textit{if }} & m+l\geq k,
\\ 
\left( P_{m}\oplus P_{l}\right) \cdot\left( \tilde{I}_{n-1}\boxplus\left( 
\tilde{I}-\left( P_{k}\oplus0\right) \right) \right) =\tilde{I}%
_{n-1}\boxplus\left( \tilde{I}-\left( P_{k-m-l}\oplus0\right) \right) & 
\text{\textit{if }} & m+l<k,%
\end{array}
\right. 
\end{equation*}
\textit{for representatives of different types and by adding up
corresponding indices }$m,l,j,n,k$ involved\textit{\ for representatives of
the same type}.

\textbf{Corollary 1}. \textit{The cancellation law holds for projections }$%
\alpha$\textit{\ over }$C\left( \mathbb{P}^{1}\left( \mathcal{T}\right)
\right) $\textit{\ of rank }$n\geq1$\textit{\ where }$n$\textit{\ is the
rank of the projection }$\sigma\left( \alpha\right) \in M_{\infty}\left(
C\left( \mathbb{T}\right) \right) $\textit{\ at any point of }$T$\textit{,
but fails for projections }$\alpha$\textit{\ over }$C\left( \mathbb{P}%
^{1}\left( \mathcal{T}\right) \right) $\textit{\ of rank} $0$.

We also get the following details about the positive cone of $K_{0}\left(
C\left( \mathbb{P}^{1}\left( \mathcal{T}\right) \right) \right) $, extending
the information provided by Corollary 3.4 of \cite{HaMaSz}.

\textbf{Corollary 2}. \textit{The positive cone of }$K_{0}\left( C\left( 
\mathbb{P}^{1}\left( \mathcal{T}\right) \right) \right) \cong\mathbb{Z}\left[
e_{11}\oplus0\right] \oplus\mathbb{Z}\left[ \tilde {I}\right] $\textit{\ is }%
\begin{equation*}
\left( \mathbb{Z}_{\geq}\left[ e_{11}\oplus0\right] \oplus0\left[ \tilde{I}%
\right] \right) \cup\left( \mathbb{Z}\left[ e_{11}\oplus0\right] \oplus%
\mathbb{Z}_{>}\left[ \tilde{I}\right] \right) . 
\end{equation*}
\textit{The canonical homomorphism from the monoid }$\mathfrak{P}\left(
C\left( \mathbb{P}^{1}\left( \mathcal{T}\right) \right) \right) $\textit{\
to }$K_{0}\left( C\left( \mathbb{P}^{1}\left( \mathcal{T}\right) \right)
\right) $\textit{\ sends}%
\begin{equation*}
\left\{ 
\begin{array}{lll}
P_{m}\oplus P_{l} & \mapsto & \left( m+l\right) \left[ e_{11}\oplus0\right]
\\ 
\tilde{I}_{n}\boxplus\left( P_{j}\oplus0\right) & \mapsto & j\left[
e_{11}\oplus0\right] \oplus n\left[ \tilde{I}\right] \\ 
\tilde{I}_{n-1}\boxplus\left( \tilde{I}-\left( P_{k}\oplus0\right) \right) & 
\mapsto & -k\left[ e_{11}\oplus0\right] \oplus n\left[ \tilde{I}\right]%
\end{array}
\right. 
\end{equation*}
\textit{for }$m,l,j\in\mathbb{Z}_{\geq}$\textit{\ and} $n,k\in\mathbb{Z}_{>}$%
.

We briefly compare the quantum complex projective space $\mathbb{P}%
^{1}\left( \mathcal{T}\right) $ with the Podle\'{s} quantum sphere $\mathbb{S%
}_{\mu c}^{2}$ for $\mu\in\left( -1,1\right) $ and $c>0$, using the groupoid
approach.

By the description of the structure of $C\left( \mathbb{S}_{\mu
c}^{2}\right) $ in \cite{Sh:qpsu}, $C\left( \mathbb{S}_{\mu c}^{2}\right) $
can be realized as the groupoid C*-algebra $C^{\ast}\left( \mathcal{F}%
\right) $ of the subgroupoid 
\begin{equation*}
\mathcal{F}:=\left\{ \left( n,n,p,q\right) :n\in\mathbb{Z}\text{ such that }%
\ \left( p,q\right) ,\left( p+n,q+n\right) \in\mathcal{G}^{\left( 0\right)
}\right\} 
\end{equation*}
of $\mathcal{G}$, sharing the same unit space $\mathcal{F}^{\left( 0\right)
}=\mathcal{G}^{\left( 0\right) }$ with $\mathcal{G}$.

The open subgroupoid 
\begin{equation*}
\mathcal{F}|_{U}=\left\{ \left( n,n,p,\infty\right) ,\left( n,n,\infty
,q\right) :n\in\mathbb{Z}\text{ such that }p,q,p+n,q+n\in\mathbb{Z}_{\geq
}\right\} 
\end{equation*}
of $\mathcal{F}$ is isomorphic to the disjoint union $\mathfrak{K}_{+}\sqcup%
\mathfrak{K}_{-}$ of two copies of the groupoid $\mathfrak{K}:=\left. \left( 
\mathbb{Z}\ltimes\mathbb{Z}\right) \right\vert _{\mathbb{Z}_{\geq}}$ under
the map $\left( n,n,p,\infty\right) \mapsto\left( n,p\right) \in\mathfrak{K}%
_{+}$ and $\left( n,n,\infty,q\right) \mapsto\left( n,q\right) \in\mathfrak{K%
}_{-}$. Thus 
\begin{equation*}
C^{\ast}\left( \mathcal{F}|_{U}\right) \cong C^{\ast}\left( \mathfrak{K}%
_{+}\right) \oplus C^{\ast}\left( \mathfrak{K}_{-}\right) \cong\mathcal{K}%
\left( \ell^{2}\left( \mathbb{Z}_{\geq}\right) \right) \oplus\mathcal{K}%
\left( \ell^{2}\left( \mathbb{Z}_{\geq}\right) \right) . 
\end{equation*}
With $\mathcal{F}^{\left( 0\right) }\backslash U=\left\{ \left(
\infty,\infty\right) \right\} $ and $\mathcal{F}|_{\left\{ \left(
\infty,\infty\right) \right\} }=\left\{ \left( n,n,\infty,\infty\right) :n\in%
\mathbb{Z}\right\} $ isomorphic to the group $\mathbb{Z}$, we get the short
exact sequence

\begin{equation*}
0\rightarrow C^{\ast}\left( \left. \mathcal{F}\right\vert _{U}\right) \cong%
\mathcal{K}\oplus\mathcal{K}\overset{\iota}{\rightarrow}C^{\ast}\left( 
\mathcal{F}\right) \cong C\left( \mathbb{S}_{\mu c}^{2}\right) \overset{%
\sigma}{\rightarrow}C^{\ast}\left( \left. \mathcal{F}\right\vert _{\left\{
\left( \infty,\infty\right) \right\} }\right) \cong C\left( \mathbb{T}%
\right) \rightarrow0, 
\end{equation*}
where $\delta_{\left( 1,1,\infty,\infty\right) }\in C^{\ast}\left( \left. 
\mathcal{F}\right\vert _{\left\{ \left( \infty,\infty\right) \right\}
}\right) $ is identified with $z:=\func{id}_{\mathbb{T}}\in C\left( \mathbb{T%
}\right) $ under $\sigma$. In the induced 6-term exact sequence of $K$%
-groups, the index homomorphism $\eta:K_{1}\left( C\left( \mathbb{T}\right)
\right) \rightarrow K_{0}\left( \mathcal{K}\oplus\mathcal{K}\right) $ sends $%
\left[ z\right] $ to $\left( -\left[ e_{11}\right] \right) \oplus\left( -%
\left[ e_{11}\right] \right) $, and hence $K_{0}\left( \iota\right) \left( %
\left[ e_{11}\oplus0\right] \right) =-K_{0}\left( \iota\right) \left( \left[
0\oplus e_{11}\right] \right) $, leading to $K_{0}\left( \iota\right) \left( %
\left[ P_{m}\oplus P_{l}\right] \right) =\left( m-l\right) \left[
e_{11}\oplus0\right] $ in $K_{0}\left( C\left( \mathbb{S}_{\mu c}^{2}\right)
\right) \cong\mathbb{Z}\left[ e_{11}\oplus0\right] \oplus\mathbb{Z}\left[ 
\tilde {I}\right] $.

By the same kind of analysis carried out above for $C\left( \mathbb{P}%
^{1}\left( \mathcal{T}\right) \right) $, we get the following results.

\textbf{Theorem 3}. \textit{The abelian monoid }$\mathfrak{P}\left( C\left( 
\mathbb{S}_{\mu c}^{2}\right) \right) $\textit{\ of unitary equivalence
classes of projections over }$C\left( \mathbb{S}_{\mu c}^{2}\right) $\textit{%
\ consists of (the representatives) }$P_{m}\oplus P_{l}$\textit{, }$\tilde{I}%
_{n}\boxplus\left( P_{j}\oplus0\right) $\textit{, and }$\tilde {I}%
_{n-1}\boxplus\left( \tilde{I}-\left( P_{k}\oplus0\right) \right) $\textit{\
(or equivalently, }$\tilde{I}_{n}\boxplus\left( 0\oplus P_{k}\right) $%
\textit{) for }$m,l,j\in\mathbb{Z}_{\geq}$\textit{\ and }$n,k\in\mathbb{Z}%
_{>}$\textit{, where }$\tilde{I}$\textit{\ is the identity of }$\left( 
\mathcal{K}\oplus\mathcal{K}\right) ^{+}\subset C\left( \mathbb{S}_{\mu
c}^{2}\right) $\textit{\ and }$P_{k}$\textit{\ is the identity of }$%
M_{k}\left( \mathbb{C}\right) \subset\mathcal{K}$\textit{, with its binary
operation }$\cdot$\textit{\ specified by}%
\begin{equation*}
\left\{ 
\begin{array}{lll}
\left( P_{m}\oplus P_{l}\right) \cdot\left( \tilde{I}_{n}\boxplus\left(
P_{j}\oplus0\right) \right) =\tilde{I}_{n}\boxplus\left(
P_{m+j-l}\oplus0\right) & \text{\textit{if }} & m+j\geq l, \\ 
\left( P_{m}\oplus P_{l}\right) \cdot\left( \tilde{I}_{n}\boxplus\left(
P_{j}\oplus0\right) \right) =\tilde{I}_{n-1}\boxplus\left( \tilde {I}-\left(
P_{l-m-j}\oplus0\right) \right) & \text{\textit{if }} & m+j<l, \\ 
\left( \tilde{I}_{n}\boxplus\left( P_{j}\oplus0\right) \right) \cdot\left( 
\tilde{I}_{n^{\prime}-1}\boxplus\left( \tilde{I}-\left( P_{k}\oplus0\right)
\right) \right) =\tilde{I}_{n+n^{\prime}}\boxplus\left( P_{j-k}\oplus0\right)
& \text{\textit{if }} & j\geq k, \\ 
\left( \tilde{I}_{n}\boxplus\left( P_{j}\oplus0\right) \right) \cdot\left( 
\tilde{I}_{n^{\prime}-1}\boxplus\left( \tilde{I}-\left( P_{k}\oplus0\right)
\right) \right) =\tilde{I}_{n+n^{\prime}-1}\boxplus\left( \tilde{I}-\left(
P_{k-j}\oplus0\right) \right) & \text{\textit{if }} & j<k, \\ 
\left( P_{m}\oplus P_{l}\right) \cdot\left( \tilde{I}_{n-1}\boxplus\left( 
\tilde{I}-\left( P_{k}\oplus0\right) \right) \right) =\tilde{I}%
_{n}\boxplus\left( P_{m-k-l}\oplus0\right) & \text{\textit{if }} & m\geq k+l,
\\ 
\left( P_{m}\oplus P_{l}\right) \cdot\left( \tilde{I}_{n-1}\boxplus\left( 
\tilde{I}-\left( P_{k}\oplus0\right) \right) \right) =\tilde{I}%
_{n-1}\boxplus\left( \tilde{I}-\left( P_{k+l-m}\oplus0\right) \right) & 
\text{\textit{if }} & m<k+l,%
\end{array}
\right. 
\end{equation*}
\textit{for representatives of different types and by adding up
corresponding indices }$m,l,j,n,k$ involved\textit{\ for representatives of
the same type}.

\textbf{Corollary 3}. \textit{The cancellation law holds for projections }$%
\alpha$\textit{\ over }$C\left( \mathbb{S}_{\mu c}^{2}\right) $\textit{\ of
rank }$n\geq1$\textit{\ where }$n$\textit{\ is the rank of the projection }$%
\sigma\left( \alpha\right) \in M_{\infty}\left( C\left( \mathbb{T}\right)
\right) $\textit{\ at any point of }$\mathbb{T}$\textit{, but fails for
projections }$\alpha$\textit{\ over }$C\left( C\left( \mathbb{S}_{\mu
c}^{2}\right) \right) $\textit{\ of rank} $0$.

The following details about the positive cone of $K_{0}\left( C\left( 
\mathbb{S}_{\mu c}^{2}\right) \right) $ extend the information provided by
Corollary 4.4 of \cite{HaMaSz}.

\textbf{Corollary 4}. \textit{The positive cone of }$K_{0}\left( C\left( 
\mathbb{S}_{\mu c}^{2}\right) \right) \cong\mathbb{Z}\left[ e_{11}\oplus0%
\right] \oplus\mathbb{Z}\left[ \tilde{I}\right] $\textit{\ is }%
\begin{equation*}
\left( \mathbb{Z}\left[ e_{11}\oplus0\right] \oplus0\left[ \tilde {I}\right]
\right) \cup\left( \mathbb{Z}\left[ e_{11}\oplus0\right] \oplus\mathbb{Z}_{>}%
\left[ \tilde{I}\right] \right) . 
\end{equation*}
\textit{The canonical homomorphism from the monoid }$\mathfrak{P}\left(
C\left( \mathbb{S}_{\mu c}^{2}\right) \right) $\textit{\ to }$K_{0}\left(
C\left( \mathbb{S}_{\mu c}^{2}\right) \right) $\textit{\ sends}%
\begin{equation*}
\left\{ 
\begin{array}{lll}
P_{m}\oplus P_{l} & \mapsto & \left( m-l\right) \left[ e_{11}\oplus0\right]
\\ 
\tilde{I}_{n}\boxplus\left( P_{j}\oplus0\right) & \mapsto & j\left[
e_{11}\oplus0\right] \oplus n\left[ \tilde{I}\right] \\ 
\tilde{I}_{n-1}\boxplus\left( \tilde{I}-\left( P_{k}\oplus0\right) \right) & 
\mapsto & -k\left[ e_{11}\oplus0\right] \oplus n\left[ \tilde{I}\right]%
\end{array}
\right. 
\end{equation*}
\textit{for }$m,l,j\in\mathbb{Z}_{\geq}$\textit{\ and} $n,k\in\mathbb{Z}_{>}$%
.

Comparing the above results, we see that quantum complex projective lines $%
\mathbb{P}^{1}\left( \mathcal{T}\right) $ and $\mathbb{S}_{\mu c}^{2}$ are
distinguished apart by the monoid structures of $\mathfrak{P}\left( C\left( 
\mathbb{P}^{1}\left( \mathcal{T}\right) \right) \right) $ and $\mathfrak{P}%
\left( C\left( \mathbb{S}_{\mu c}^{2}\right) \right) $, and also by the
positive cone of their $K_{0}$-groups.

\section{Line bundles over quantum projective line}

In this section, we identify the quantum line bundles $L_{k}:=C\left( 
\mathbb{S}_{H}^{3}\right) _{k}$ of degree $k$ over $C\left( \mathbb{P}%
^{1}\left( \mathcal{T}\right) \right) $ with a concrete (unitary equivalence
class of) projection classified in the previous section.

To distinguish between ordinary function product and convolution product, we
denote the groupoid C*-algebraic multiplication of elements in $C^{\ast
}\left( \mathcal{G}\right) \supset C_{c}\left( \mathcal{G}\right) $ by $\ast$%
, while omitting $\ast$ when the elements are presented as operators or when
they are multiplied together pointwise as functions.

Recall from earlier section that $L_{k}=\overline{C_{c}\left( \mathcal{G}%
_{k}\right) }\subset C\left( \mathbb{S}_{H}^{3}\right) $ where 
\begin{equation*}
\mathcal{G}_{k}:=\left\{ \left( n+k,-n,p,q\right) :n\in\mathbb{Z}\text{ such
that }\ \left( p,q\right) ,\left( p+n+k,q-n\right) \in\mathcal{G}^{\left(
0\right) }\right\} \subset\mathcal{G}, 
\end{equation*}
and $L_{0}=\overline{C_{c}\left( \mathcal{G}_{0}\right) }=C\left( \mathbb{P}%
^{1}\left( \mathcal{T}\right) \right) $. Furthermore the groupoid C*-algebra 
$C^{\ast}\left( \mathcal{G}\right) \cong C\left( \mathbb{S}_{H}^{3}\right) =%
\overline{\oplus_{k\in\mathbb{Z}}L_{k}}$ is a $\mathbb{Z}$-graded algebra.

Let $k>0$. We identify below separately $L_{k}$ and $L_{-k}$ with a
representative of projections over $C\left( \mathbb{P}^{1}\left( \mathcal{T}%
\right) \right) $ classified earlier.

The characteristic function $\chi_{A}\in C_{c}\left( \mathcal{G}_{k}\right)
\subset C\left( \mathbb{S}_{H}^{3}\right) $ of the compact set 
\begin{equation*}
A:=\left\{ \left( k,0,p,q\right) :\left( p,q\right) ,\left( p+k,q\right) \in%
\mathcal{G}^{\left( 0\right) }\right\} \subset\mathcal{G}
\end{equation*}
is an isometry with $\chi_{A}^{\ast}\ast\chi_{A}=\chi_{\mathcal{G}^{\left(
0\right) }}\equiv1\in C^{\ast}\left( \mathcal{G}_{0}\right) $ and 
\begin{equation*}
\chi_{A}\ast\chi_{A}^{\ast}=\chi_{B}=\tilde{I}-\left( P_{k}\oplus0\right) 
\end{equation*}
a projection in $C^{\ast}\left( \mathcal{G}_{0}\right) $ $\equiv C\left( 
\mathbb{P}^{1}\left( \mathcal{T}\right) \right) $for the set 
\begin{equation*}
B:=\left\{ \left( 0,0,p^{\prime},q^{\prime}\right) \in\mathcal{G}:\
p^{\prime}\geq k,\ q^{\prime}\geq0\right\} . 
\end{equation*}
So with $\chi_{B}\ast\chi_{A}=\chi_{A}\in C_{c}\left( \mathcal{G}_{k}\right) 
$, we get a left $C^{\ast}\left( \mathcal{G}_{0}\right) $-module homomorphism%
\begin{equation*}
x\in C^{\ast}\left( \mathcal{G}_{0}\right) \ast\chi_{B}\mapsto x\ast\chi
_{A}\in\overline{C_{c}\left( \mathcal{G}_{k}\right) }\equiv L_{k}
\end{equation*}
with well-defined inverse map 
\begin{equation*}
y\in\overline{C_{c}\left( \mathcal{G}_{k}\right) }\mapsto y\ast\chi
_{A}^{\ast}=y\ast\chi_{A}^{\ast}\ast\chi_{B}\in C^{\ast}\left( \mathcal{G}%
_{0}\right) \ast\chi_{B}
\end{equation*}
since $\chi_{A}^{\ast}\in C_{c}\left( \mathcal{G}_{-k}\right) $ and hence $%
C_{c}\left( \mathcal{G}_{k}\right) \ast\chi_{A}^{\ast}\subset C_{c}\left( 
\mathcal{G}_{0}\right) $. Now $L_{k}$ being isomorphic to the left $C^{\ast
}\left( \mathcal{G}_{0}\right) $-module $C^{\ast}\left( \mathcal{G}%
_{0}\right) \left( \tilde{I}-\left( P_{k}\oplus0\right) \right) $ is
identified with the rank-one projection $\tilde{I}-\left(
P_{k}\oplus0\right) $.

Next we show that in the left $C^{\ast}\left( \mathcal{G}_{0}\right) $%
-module decomposition%
\begin{equation*}
L_{-k}=L_{-k}\ast\chi_{B}\oplus L_{-k}\ast\left( 1-\chi_{B}\right) \equiv
L_{-k}\left( \tilde{I}-\left( P_{k}\oplus0\right) \right) \oplus
L_{-k}\left( P_{k}\oplus0\right) 
\end{equation*}
by the projection $\chi_{B}$, the first component $L_{-k}\ast\chi_{B}$ can
be identified with the projection $\tilde{I}$ and the second component $%
L_{-k}\ast\left( 1-\chi_{B}\right) $ can be identified with the projection $%
P_{k}\oplus0$, leading to the conclusion that $L_{-k}$ can be identified
with the projection $\tilde{I}\boxplus\left( P_{k}\oplus0\right) $.

Indeed since $\chi_{A}\ast\chi_{A}^{\ast}=\chi_{B}$ and $\chi_{A}^{\ast}\ast%
\chi_{A}=1\equiv\chi_{\mathcal{G}^{\left( 0\right) }}$ with $\chi
_{A}^{\ast}\in C_{c}\left( \mathcal{G}_{-k}\right) $, the map%
\begin{equation*}
x\in C^{\ast}\left( \mathcal{G}_{0}\right) \mapsto x\ast\chi_{A}^{\ast
}=x\ast\chi_{A}^{\ast}\ast\chi_{B}\in\overline{C_{c}\left( \mathcal{G}%
_{-k}\right) }\ast\chi_{B}\equiv L_{-k}\ast\chi_{B}
\end{equation*}
is a $C^{\ast}\left( \mathcal{G}_{0}\right) $-module isomorphism with
inverse $y\mapsto y\ast\chi_{A}$ and hence $L_{-k}\ast\chi_{B}$ is
identified with the projection $\tilde{I}$.

On the other hand, comparing 
\begin{align*}
C_{c}\left( \mathcal{G}_{-k}\right) \ast\left( 1-\chi_{B}\right) &
=C_{c}\left( \left\{ \left( n-k,-n,p,\infty\right) :0\leq p<k\text{ and}\
p+n\geq k\right\} \right) \\
& =C_{c}\left( \left\{ \left( n,-n-k,p,\infty\right) :0\leq p<k\text{ and}\
p+n\geq0\right\} \right)
\end{align*}
and%
\begin{equation*}
C_{c}\left( \mathcal{G}_{0}\right) \ast\left( 1-\chi_{B}\right) =C_{c}\left(
\left\{ \left( n,-n,p,\infty\right) :0\leq p<k\text{\ and}\ p+n\geq0\right\}
\right) 
\end{equation*}
where with the last coordinate being $\infty$, the second coordinate becomes
irrelevant, we get a $\left( C_{c}\left( \mathcal{G}_{0}\right) ,\ast\right) 
$-module isomorphism%
\begin{equation*}
f\in C_{c}\left( \mathcal{G}_{0}\right) \ast\left( 1-\chi_{B}\right) \mapsto
f\circ\tau\in C_{c}\left( \mathcal{G}_{-k}\right) \ast\left(
1-\chi_{B}\right) 
\end{equation*}
where $\tau\left( n,-n-k,p,\infty\right) :=\left( n,-n,p,\infty\right) $,
which extends to a $C^{\ast}\left( \mathcal{G}_{0}\right) $-module
isomorphism 
\begin{equation*}
C^{\ast}\left( \mathcal{G}_{0}\right) \ast\left( 1-\chi_{B}\right) \equiv
C^{\ast}\left( \mathcal{G}_{0}\right) \left( P_{k}\oplus0\right) \rightarrow%
\overline{C_{c}\left( \mathcal{G}_{-k}\right) }\ast\left( 1-\chi_{B}\right)
\equiv L_{-k}\ast\left( 1-\chi_{B}\right) . 
\end{equation*}
So the $C^{\ast}\left( \mathcal{G}_{0}\right) $-module $L_{-k}\ast\left(
1-\chi_{B}\right) $ is identified with the projection $P_{k}\oplus0$.

We summarize as follows.

\textbf{Theorem 4}. \textit{The quantum line bundle }$L_{k}\equiv C\left( 
\mathbb{S}_{H}^{3}\right) _{k}$\textit{\ of degree }$k\in\mathbb{Z}$\textit{%
\ over }$C\left( \mathbb{P}^{1}\left( \mathcal{T}\right) \right) $\textit{\
is isomorphic to the finitely generated projective left module over }$%
C\left( \mathbb{P}^{1}\left( \mathcal{T}\right) \right) $\textit{\
determined by the projection }$\tilde{I}-\left( P_{k}\oplus0\right) $\textit{%
\ if }$k\geq 0$\textit{, and the projection }$\tilde{I}\boxplus\left(
P_{-k}\oplus0\right) $\textit{\ if} $k<0$.

\textbf{Corollary 5}. \textit{The quantum line bundles }$L_{k}\equiv C\left( 
\mathbb{S}_{H}^{3}\right) _{k}$\textit{\ with }$k\in\mathbb{Z}$\textit{\
provide a complete list of mutually non-isomorphic rank-one finitely
generated left projective modules over} $C\left( \mathbb{P}^{1}\left( 
\mathcal{T}\right) \right) $.

It is interesting to note that in the case of quantum teardrops $%
WP_{q}\left( k,l\right) $, the quantum principal $U\left( 1\right) $-bundles 
$\mathcal{L}\left( k\right) $ of degree $k$ over $C\left( WP_{q}\left(
k,l\right) \right) $ introduced by Brzezi\'{n}ski and Fairfax \cite{BrzeFair}
do not exhaust all rank-one finitely generated projective modules over $%
C\left( WP_{q}\left( k,l\right) \right) $ by the result of \cite{Sh:qt}.

\end{document}